\documentclass{amsart}
\usepackage{lastpage}
\usepackage{yhmath,verbatim,mathabx}
\usepackage[all]{xy}

\DeclareFontFamily{U}{mathx}{\hyphenchar\font45}
\DeclareFontShape{U}{mathx}{m}{n}{
	<5> <6> <7> <8> <9> <10>
	<10.95> <12> <14.4> <17.28> <20.74> <24.88>
	mathx10
}{}
\DeclareSymbolFont{mathx}{U}{mathx}{m}{n}
\DeclareFontSubstitution{U}{mathx}{m}{n}
\DeclareMathAccent{\widecheck}{0}{mathx}{"71}
\DeclareMathAccent{\wideparen}{0}{mathx}{"75}

\newcommand{\bdis}{\begin{displaymath}}
\newcommand{\edis}{\end{displaymath}}
\newcommand{\be}{\begin{equation}}
\newcommand{\ee}{\end{equation}}
\newcommand{\mbb}{\mathbb}
\newcommand{\mcal}{\mathcal}

\newcommand{\zf}{\zeta\left(\frac{1}{2}+it\right)}

\DeclareMathOperator{\dn}{dn} 

\newcommand{\Opar}{\left\{1+\mathcal{O}\left(\frac{\ln\ln L}{\ln L}\right)\right\}}

\DeclareMathOperator{\sn}{sn}
\DeclareMathOperator{\cn}{cn}

\theoremstyle{definition}

\theoremstyle{remark}
\newtheorem{remark}[]{Remark}

\newtheorem*{mydef1}{{\bf Theorem}}

\newtheorem*{mydef51}{{\bf Lemma 1}}

\newtheorem*{mydef52}{{\bf Lemma 2}}

\newtheorem*{mydef53}{{\bf Lemma 3}}

\newtheorem*{mydef54}{{\bf Lemma 4}}

\newtheorem*{mydef55}{{\bf Lemma 5}}

\numberwithin{equation}{section}

\begin{document}

\title{Jacob's ladders and completely new exact synergetic formula for Jacobi's elliptic functions together with Bessel's functions excited by the function $|\zf|^2$} 

\author{Jan Moser}

\address{Department of Mathematical Analysis and Numerical Mathematics, Comenius University, Mlynska Dolina M105, 842 48 Bratislava, SLOVAKIA}

\email{jan.mozer@fmph.uniba.sk}


\begin{abstract}
In this paper we obtain a set of five new transmutations of the mother formula. Further, we obtain the second set of ten exact metafunctional equations by crossbreeding on every two elements of the previous set. Elements of the last set represent a kind of synergetic formulae describing cooperative interactions between corresponding sets of values of basic classical functions.  

\end{abstract}
\maketitle

\section{Introduction} 

\subsection{} 

Let us remind that the following set of values 
\be \label{1.1} 
\begin{split}
& \left\{\left|\zf\right|^2\right\}, \\ 
& \{f_1(t)\}=\{\sin^2t\},\ \{f_2(t)\}=\{\cos^2t\},\ \{f_3(t)\}=\{\cos 2t\}, \\ 
& t\in [\pi L,\pi L+U],\ U\in(0,\pi/4),\ L\in\mbb{N} 
\end{split}
\ee 
generate the following secondary asymptotic complete hybrid formula (see \cite{7}, (3.7), \cite{6} -- \cite{9},  $k_1=k_2=k_3=1$) 
\be \label{1.2} 
\begin{split}
& \left|\zeta\left(\frac 12+i\alpha_1^{1,1}\right)\right|^2\sin^2\alpha_0^{1,1}-\\ 
& - \left\{1+\mcal{O}\left(\frac{\ln\ln L}{\ln L}\right)\right\}\left|\zeta\left(\frac 12+i\alpha_1^{2,1}\right)\right|^2\cos^2\alpha_0^{2,1}+\\ 
& + \left|\zeta\left(\frac 12+i\alpha_1^{3,1}\right)\right|^2\cos(2\alpha_0^{3,1})=0, \\ 
& \forall\- L\leq L_0>0, 
\end{split}
\ee 
($L_0$ is sufficiently big), where 
\be \label{1.3} 
\begin{split}
& \alpha_r^{l,1}=\alpha_r(U,\pi L,f_l),\ r=0,1,\ l=1,2,3, \\ 
& \alpha_0^{l,1}\in (\pi L,\pi L+U),\ \alpha_1^{l,1}\in (\overset{1}{\wideparen{\pi L}},\overset{1}{\wideparen{\pi L+U}}), 
\end{split}
\ee 
and 
\bdis 
[\overset{1}{\wideparen{\pi L}},\overset{1}{\wideparen{\pi L+U}}]
\edis 
is the first reverse iteration of the basic segment 
\bdis 
[\pi L,\pi L+U]=[\overset{0}{\wideparen{\pi L}},\overset{0}{\wideparen{\pi L+U}}]
\edis 
by means of the Jacob's ladder, see \cite{2}, \cite{3}. 

\subsection{} 

Next, we have obtained in our paper \cite{8} the set of six transmutations of the mother formula (\ref{1.2}). For example, there are continuum sets 
\be \label{1.4} 
\overset{(1)}{\Omega}_l(U,L;[|\zeta_{0,5}|^2]),\ \overset{(7)}{\Omega}_l(U,L;[|\zeta_{0,5}|^2])\subset\mbb{C},\ l=1,2,3 
\ee  
such that for every of elements 
\bdis 
s_l^1\in \overset{(1)}{\Omega}_l,\ s_l^7\in\overset{(7)}{\Omega}_l 
\edis  
we have the following transmutation of (\ref{1.2}) 
\be \label{1.5} 
\begin{split}
& |\zeta(s_1^1)|^2|\sn(s_1^7,k_1)|^2-\\ 
& -\{1+\mcal{O}(\frac{\ln\ln L}{\ln L})\}|\zeta(s_2^1)|^2|\cn(s_2^7,k_2)|^2+\\ 
& + |\zeta(s_3^1)|^2|\dn(s_3^7,k_3)|=0, \\ 
& [(k_1)^2,(k_2)^2,(k_3)^2]\in [0,1]^3, 
\end{split}
\ee 
(comp. also the connection between the basic functions in \cite{9}, (1.6) -- (1.11)).  

\subsection{} 

Let us remind that the factor 
\be \label{1.6} 
1+\mcal{O}\left(\frac{\ln\ln L}{\ln L}\right) 
\ee  
which is contained in every element of the set of seven transmutations is the identical one. Consequently, we have used in our paper \cite{9} the operation of crossbreeding (see \cite{4} -- \cite{7}, i.e. in our case the elimination of the function (\ref{1.6})) on every two different elements of the set of transmutations. We have obtained as result the first set of fifteen exact meta-functional equations. For example, (see \cite{9}, (3.8)): 
\be \label{1.7} 
\begin{split}
& |\zeta(s_1^1)|^2|\zeta(s_1^2)|^2|\cn(s_2^7,k_2)|^2+|\zeta(s_3^1)|^2|\zeta(s_3^2)|^2|\cn(s_2^7,k_2)|^2= \\ 
& |\zeta(s_1^1)|^2|\zeta(s_2^2)|^2|\sn(s_1^7,k_1)|^2+|\zeta(s_3^1)|^2|\zeta(s_3^2)|^2|\dn(s_3^7,k_3)|. 
\end{split}
\ee  

\subsection{} 

In this paper we obtain a set of five new transmutations of the mother formula. Further, by crossbreeding on every two elements of this set we obtain the second set of ten exact meta-functional equations. For example, there are the sets 
\be \label{1.8} 
\begin{split}
& \overset{(6)}{\Omega}_l(U,L,p_l;[|\zeta_{0,5}|^2]),\ \overset{(7)}{\Omega}_l(U,L,k_l;[|\zeta_{0,5}|^2]), \\
& \overset{(11)}{\Omega}_l(U,L,p_{l+3};[|\zeta_{0,5}|^2]),\ \overset{(12)}{\Omega}_l(U,L,k_{l+3};[|\zeta_{0,5}|^2])\subset\mbb{C}, \\ 
& l=1,2,3 
\end{split}
\ee 
such that for every 
\be \label{1.9} 
s_l^n=s_l^n([|\zeta_{0,5}|^2])\in \overset{(n)}{\Omega}_l,\ 
l=1,2,3,\ n=6,7,11,12 
\ee  
we have the following exact meta-functional equation: 
\be \label{1.10} 
\begin{split}
& |J_{p_1}(s_1^6)|^2|\cn(s_2^7,k_2)|^2|J_{p_4}(s_1^{11})|^2|\cn(s_2^{12},k_5)|^2+ \\ 
& |J_{p_3}(s_3^6)|^2|\cn(s_2^7,k_2)|^2|J_{p_6}(s_3^{11})|^2|\cn(s_2^{12},k_5)|^2=\\ 
& = |J_{p_2}(s_2^6)|^2|\sn(s_1^7,k_1)|^2|J_{p_2}(s_2^{11})|^2|\sn(s_1^{12},k_4)|^2+ \\ 
& |J_{p_2}(s_2^6)|^2|\dn(s_3^7,k_3)||J_{p_5}(s_2^{11})|^2|\dn(s_3^{12},k_6)|^2, \\ 
& (p_1,p_2,p_3), (p_4,p_5,p_6)\in\mbb{Z}^3,\\ 
& [(k_1)^2,(k_2)^2,(k_3)^2],[(k_4)^2,(k_5)^2,(k_6)^2]\in (0,1)^3. 
\end{split}
\ee 

\begin{remark}
Though the formula (\ref{1.10}) does not contain the external factors of type (comp. (\ref{1.7})) 
\bdis 
|\zeta(s_1^1)|^2,\ |\zeta(s_2^1)|^2, \dots ,
\edis 
for example, even so this one is generated (internally) mainly by the function 
\bdis 
\left|\zf\right|^2=|\zeta_{0,5}|^2, 
\edis  
(comp. (\ref{1.8}), (\ref{1.9}) and \cite{9}, (1.6) -- (1.11)). 
\end{remark} 

\begin{remark}
The morphology (external structure for us) of the formula (\ref{1.10}), as well as for all other formulae in the Theorem, is unambiguous determined by the mother formula (\ref{1.2}) and by corresponding operation of crossbreeding. 
\end{remark} 

\begin{remark}
Consequently, the following is true: all ten exact meta-functional equations obtained in this paper are the direct descendants of the mother formula (\ref{1.2}). 
\end{remark} 

\begin{remark}
From the point of view of our $\zeta$-alchemy, the global synergetic formula (\ref{1.10}) represents the result ($\zeta$-compound ) of cooperative interactions excited by the function 
\bdis 
\left|\zf\right| 
\edis 
between the sets (substances) 
\be \label{1.11} 
\begin{split}
& \{|J_{p_l}(s)|\},\ \{|J_{p_{l+3}}(s)|\},\ l=1,2,3, \\ 
& \{|\sn(s,k_1)|\},\ \{|\cn(s,k_2)|\},\ \{|\dn(s,k_3)|\}, \\ 
& \{|\sn(s,k_4)|\},\ \{|\cn(s,k_5)|\},\ \{|\dn(s,k_6)|\}, \\ 
& s\in\mbb{C}\setminus \{P\}, 
\end{split}
\ee  
where $\{P\}$ is the set of corresponding poles. 
Of course, we may write another nine classes on basis of our Theorem in Section 3 of this paper. 
\end{remark}

\begin{remark}
The synergetic formula (\ref{1.10}) represents the completely new type of formula simultaneously for the following theories: 
\begin{itemize}
	\item[(a)] Riemann's zeta-function, 
	\item[(b)] Jacobi's elliptic functions, 
	\item[(c)] Bessel's functions. 
\end{itemize}
\end{remark} 

\begin{remark}
This paper is also based on new notions and methods in the theory of the Riemann's function we have introduced in our series of 51 papers concerning Jacob's ladders. These can be found in arXiv[math.CA] starting with the paper \cite{1}. 
\end{remark} 

\section{Further infinite set of transmutation of the mother formula} 

\subsection{} 

We define the level curves (comp. \cite{9}, (1.10), (1.11)) 
\be \label{2.1} 
\overset{(8)}{\Omega}_l=\overset{(8)}{\Omega}_l(U,L;[|\zeta_{0,5}|^2]),\ l=1,2,3 
\ee  
as the loci 
\be \label{2.2} 
\begin{split}
& |\cos s_1^8|=\left|\zeta\left(\frac 12+i\alpha_1^{1,1}\right)\right|, \\ 
& |\cos s_2^8|=\left|\zeta\left(\frac 12+i\alpha_1^{2,1}\right)\right|, \\ 
& |\cos s_3^8|=\left|\zeta\left(\frac 12+i\alpha_1^{3,1}\right)\right|; \\ 
& s_l^8\in \overset{(8)}{\Omega}_l, 
\end{split}
\ee 
for every fixed and admissible $U,L$, correspondingly. Now we have, see (\ref{1.2}), (\ref{2.2}) and \cite{8}, (4.3), the following 

\begin{mydef51}
For every fixed and admissible $U,L$ there are sets 
\be \label{2.3} 
\overset{(3)}{\Omega}_l(U,L),\ \overset{(8)}{\Omega}_l(U,L)\subset\mbb{C},\ l=1,2,3 
\ee  
such that the following formula holds true (transmutation of (\ref{1.2})): 
\be \label{2.4} 
\begin{split}
& |\cos s_1^8|^2|\cos s_1^3|^2-\left\{1+\mcal{O}\left(\frac{\ln\ln L}{\ln L}\right)\right\}|\cos s_2^3|^2|\cos s_2^8|^2+ \\ 
& + |\cos s_3^8|^2|\cos s_3^3|=0. 
\end{split}
\ee 
\end{mydef51} 

\subsection{} 

Let us consider now the functions 
\be \label{2.5} 
(s)^{n_4},\ (s)^{n_5},\ (s)^{n_6},\ \forall\- (n_4,n_5,n_6)\in\mbb{N}^3. 
\ee 
We define for them the level curves 
\be \label{2.6} 
\overset{(9)}{\Omega}_l=\overset{(9)}{\Omega}_l(U,L,n_{l+3};[|\zeta_{0,5}|^2]),\ l=1,2,3 
\ee  
as the loci 
\be \label{2.7} 
\begin{split}
& |s_1^9|^{n_4}=\left|\zeta\left(\frac 12+i\alpha_1^{1,1}\right)\right|, \\ 
& |s_2^9|^{n_5}=\left|\zeta\left(\frac 12+i\alpha_1^{2,1}\right)\right|, \\ 
& |s_3^9|^{n_6}=\left|\zeta\left(\frac 12+i\alpha_1^{3,1}\right)\right|, 
\end{split}
\ee 
for every fixed and admissible $U,L,n_4,n_5,n_6$, correspondingly. Now we have, see (\ref{1.2}), (\ref{2.7}) and \cite{8}, (4.8), the following 

\begin{mydef52}
For every fixed and admissible $U,L$ and for every fixed 
\bdis 
(n_1,n_2,n_3),\ (n_4,n_5,n_6)\in\mbb{N}^3 
\edis 
there are the sets 
\be \label{2.8} 
\overset{(4)}{\Omega}_l(U,L,n_l),\ \overset{(9)}{\Omega}_l(U,L,n_{l+3})\subset\mbb{C},\ l=1,2,3 
\ee  
such that we have the following formula (transmutation of (\ref{1.2})): 
\be \label{2.9} 
\begin{split}
& |s_1^4|^{2n_1}|s_1^9|^{2n_4}-\Opar |s_2^4|^{2n_2}|s_2^9|^{2n_5}+\\ 
& + |s_3^9|^{2n_6}|s_3^4|^{n_3}=0. 
\end{split}
\ee 
\end{mydef52} 

\subsection{} 

We define the level curves 
\be \label{2.10} 
\overset{(10)}{\Omega}_l=\overset{(8)}{\Omega}_l(U,L;[|\zeta_{0,5}|^2]),\ l=1,2,3  
\ee 
as the loci 
\be \label{2.11} 
\begin{split}
& \frac{1}{|\Gamma(s_1^{10})|}=\left|\zeta\left(\frac 12+i\alpha_1^{1,1}\right)\right|, \\ 
& \frac{1}{|\Gamma(s_2^{10})|}=\left|\zeta\left(\frac 12+i\alpha_1^{2,1}\right)\right|, \\ 
& \frac{1}{|\Gamma(s_3^{10})|}=\left|\zeta\left(\frac 12+i\alpha_1^{3,1}\right)\right|, 
\end{split}
\ee  
for every fixed and admissible $U,L$, correspondingly. Now we have (see (\ref{1.2}), (\ref{2.11}) and \cite{8}, (5.3)) the following: 

\begin{mydef53}
For every fixed and admissible $U,L$ there are sets 
\be \label{2.12} 
\overset{(5)}{\Omega}_l(U,L), \overset{(10)}{\Omega}_l(U,L)\subset\mbb{C},\ l=1,2,3 
\ee 
such that the following formula (transmutation of (\ref{1.2})) holds true: 
\be \label{2.13} 
\begin{split}
& \frac{1}{|\Gamma(s_1^5)\Gamma(s_1^{10})|^2}-\Opar \frac{1}{|\Gamma(s_2^5)\Gamma(s_2^{10})|^2}+ \\ 
& + \frac{1}{|\Gamma(s_3^5)\Gamma(s_3^{10})|}=0. 
\end{split}
\ee 
\end{mydef53} 

\subsection{} 

We define, for Bessel's functions 
\be \label{2.14} 
J_{p_4}(s),\ J_{p_5}(s),\ J_{p_6}(s),\ s\in\mbb{C},\ (p_4,p_5,p_6)\in\mbb{Z}^3, 
\ee 
the level curves 
\be \label{2.15} 
\overset{(11)}{\Omega}_l=\overset{(11)}{\Omega}_l(U,L,p_{l+3};[|\zeta_{0,5}|^2]),\ l=1,2,3  
\ee 
as the loci 
\be \label{2.16} 
\begin{split}
& |J_{p_4}(s_1^{11})|=\left|\zeta\left(\frac 12+i\alpha_1^{1,1}\right)\right|, \\ 
& |J_{p_5}(s_2^{11})|=\left|\zeta\left(\frac 12+i\alpha_1^{2,1}\right)\right|, \\ 
& |J_{p_6}(s_3^{11})|=\left|\zeta\left(\frac 12+i\alpha_1^{3,1}\right)\right| 
\end{split}
\ee 
for fixed and admissible $U,L,p_4,p_5,p_6$, correspondingly. Now we have (see (\ref{1.2}), (\ref{2.16}) and \cite{8}, (5.8)) the following: 

\begin{mydef54}
For every fixed and admissible $U,L$ and for every fixed 
\bdis 
(p_1,p_2,p_3),\ (p_4,p_5,p_6)\in\mbb{Z}^3 
\edis  
there are sets 
\be \label{2.17} 
\overset{(6)}{\Omega}_l(U,L,p_l),\ \overset{(11)}{\Omega}_l(U,L,p_{l+3})\subset\mbb{C},\ l=1,2,3 
\ee  
such that the following formula (transmutation of (\ref{1.2})) holds true: 
\be \label{2.18} 
\begin{split}
& |J_{p_1}(s_1^6)|^2|J_{p_4}(s_1^{11})|^2-\Opar |J_{p_2}(s_2^6)|^2|J_{p_5}(s_2^{11})|^2+ \\ 
& + |J_{p_3}(s_3^6)||J_{p_6}(s_3^{11})|=0. 
\end{split}
\ee 
\end{mydef54} 

\subsection{} 

Fot the Jacobi's elliptic functions 
\be \label{2.19} 
\begin{split}
& \sn(s,k_4),\ \cn(s,k_5),\ \dn(s,k_6), \\ 
& [(k_4)^2,(k_5)^2,(k_6)^2]\in (0,1)^3, 
\end{split}
\ee 
we define the level curves 
\be \label{2.20} 
\overset{(12)}{\Omega}_l=\overset{(12)}{\Omega}_l(U,L,k_{l+3};[|\zeta_{0,5}|^2]),\ l=1,2,3  
\ee  
as the loci 
\be \label{2.21} 
\begin{split}
	& |\sn(s_1^{12},k_4)|=\left|\zeta\left(\frac 12+i\alpha_1^{1,1}\right)\right|, \\ 
	& |\cn(s_2^{12},k_5)|=\left|\zeta\left(\frac 12+i\alpha_1^{2,1}\right)\right|, \\ 
	& |\dn(s_3^{12},k_6)|=\left|\zeta\left(\frac 12+i\alpha_1^{3,1}\right)\right| 
\end{split}
\ee  
for fixed and admissible $U,L,k_4,k_5,k_6$, correspondingly. Now we have (see (\ref{1.2}), (\ref{2.21}) and \cite{8}, (5.13)) the following:  

\begin{mydef55}
For every fixed and admissible $U,L$ and for every fixed 
\bdis 
[(k_1)^2,(k_2)^2,(k_3)^2],\ [(k_4)^2,(k_5)^2,(k_6)^2]\in(0,1)^3 
\edis  
there are sets 
\be \label{2.22} 
\overset{(7)}{\Omega}_l(U,L,k_l),\ \overset{(12)}{\Omega}_l(U,L,k_{l+3})\subset\mbb{C},\ l=1,2,3 
\ee 
such that the following formula (transmutation of (\ref{1.2})) holds true: 
\be \label{2.23} 
\begin{split}
& |\sn(s_1^7,k_1)|^2|\sn(s_1^{12},k_4)|^2- \\ 
& - \Opar |\cn(s_2^7,k_2)|^2|\cn(s_2^{12},k_5)|^2+ \\ 
& + |\dn(s_3^7,k_3)||\dn(s_3^{12},k_6)|=0. 
\end{split}
\ee 
\end{mydef55} 

\section{List of exact meta-functional equations of the second generation} 

We have obtained the following set 
\be \label{3.1} 
\{(2.4), (2.9), (2.13),(2.18), (2.23)\}
\ee 
of five further transmutations of the mother formula (\ref{1.2}). Let us remind that the factor 
\be \label{3.2} 
\Opar 
\ee  
(comp. \cite{9}, (2.21), (2.22)) which is contained in every element of the set (\ref{3.1}) is the identical one. Consequently, we may apply the operation of crossbreeding (see \cite{4} -- \cite{7}, here the elimination of the function (\ref{3.2})) on every two different elements of the set (\ref{3.1}). 

\begin{remark}
Let the symbol 
\bdis 
(2.4)\times (2.9) \Rightarrow 
\edis 
stand for \emph{we obtain by crossbreeding of the transmutations of (2.4) and (2.9)}. 
\end{remark} 

We obtain the following Theorem as the result of crossbreeding on the set (\ref{3.1}):  

\begin{mydef1}
There are sets (see \cite{9}, (1.11)) 
\bdis 
\begin{split}
& \overset{(3)}{\Omega}_l(U,L;[|\zeta_{0,5}|^2]),\ \overset{(4)}{\Omega}_l(U,L,n_l;[|\zeta_{0,5}|^2]), \\ 
& \overset{(5)}{\Omega}_l(U,L;[|\zeta_{0,5}|^2]),\ \overset{(6)}{\Omega}_l(U,L,p_l;[|\zeta_{0,5}|^2]), \\ 
& \overset{(7)}{\Omega}_l(U,L,k_l;[|\zeta_{0,5}|^2])
\end{split}
\edis  
and, further, (see (\ref{2.1}), (\ref{2.6}), (\ref{2.10}), (\ref{2.15}) and (\ref{2.20})) 
\be \label{3.3} 
\begin{split}
	& \overset{(8)}{\Omega}_l(U,L;[|\zeta_{0,5}|^2]),\ \overset{(9)}{\Omega}_l(U,L,n_{l+3};[|\zeta_{0,5}|^2]), \\ 
	& \overset{(10)}{\Omega}_l(U,L;[|\zeta_{0,5}|^2]),\ \overset{(11)}{\Omega}_l(U,L,p_{l+3};[|\zeta_{0,5}|^2]), \\ 
	& \overset{(12)}{\Omega}_l(U,L,k_{l+3};[|\zeta_{0,5}|^2])\subset\mbb{C},\ l=1,2,3 
\end{split}
\ee  
where 
\bdis 
\begin{split}
& U\in (0,\pi/4),\ L\geq L_0>0, \\ 
& (n_1,n_2,n_3), (n_4,n_5,n_6)\in\mbb{N}^3, \\ 
& (p_1,p_2,p_3), (p_4,p_5,p_6)\in\mbb{Z}^3,\\ 
& [(k_1)^2,(k_2)^2,(k_3)^2], [(k_4)^2,(k_5)^2,(k_6)^2]\in (0,1)^3, 
\end{split}
\edis  
such that for each of the elements 
\bdis 
s_l^n\in \overset{(n)}{\Omega}_l,\ l=1,2,3,\ n=3,4,\dots,12 
\edis  
we have the following set of ten exact meta-functional equations: 
\be \label{3.4} 
\begin{split}
& (2.4)\times (2.9) \Rightarrow \\ 
& |\cos s_1^3\cos s_1^8|^2|s_2^4|^{2n_2}|s_2^9|^{2n_6}+ \\ 
& + |\cos s_3^3||\cos s_3^8|^2|s_2^4|^{2n_2}|s_2^9|^{2n_6}= \\ 
& |\cos s_2^3\cos s_2^8|^2|s_1^4|^{2n_1}|s_1^9|^{2n_4}+ \\ 
& + |\cos s_2^3\cos s_2^8|^2|s_3^4|^{n_3}|s_3^9|^{2n_6}, 
\end{split}
\ee  

\be \label{3.5} 
\begin{split}
& (2.4)\times (2.13) \Rightarrow \\ 
& \frac{|\cos s_1^3\cos s_1^8|^2}{|\Gamma(s_2^5)\Gamma(s_2^{10})|^2}+ 
\frac{|\cos s_3^3||\cos s_3^8|^2}{|\Gamma(s_2^5)\Gamma(s_2^{10})|^2}= \\ 
& =  \frac{|\cos s_2^3\cos s_2^8|^2}{|\Gamma(s_1^5)\Gamma(s_1^{10})|^2}+ 
\frac{|\cos s_2^3\cos s_2^8|^2}{|\Gamma(s_3^5)||\Gamma(s_3^{10})|^2}, 
\end{split}
\ee  

\be \label{3.6} 
\begin{split}
& (2.4)\times (2.18) \Rightarrow \\  
& |J_{p_2}(s_2^6)J_{p_5}(s_2^{11})|^2|\cos s_1^3\cos s_1^8|^2+ \\ 
& |J_{p_3}(s_2^6)J_{p_5}(s_2^{11})|^2|\cos s_3^3||\cos s_3^8|^2= \\ 
& = |J_{p_1}(s_1^6)J_{p_4}(s_1^{11})|^2|\cos s_2^3\cos s_2^8|^2+ \\ 
& + |J_{p_3}(s_3^6)||J_{p_6}(s_3^{11})|^2|\cos s_2^3\cos s_2^8|^2, 
\end{split}
\ee  

\be \label{3.7} 
\begin{split}
& (2.4)\times (2.23) \Rightarrow \\  
& |\cn(s_2^7,k_2)\cn(s_2^{12},k_5)|^2|\cos s_1^3\cos s_1^8|^2+ \\ 
& + |\cn(s_2^7,k_2)\cn(s_2^{12},k_5)|^2|\cos s_3^3||\cos s_3^8|^2= \\ 
& = |\sn(s_1^7,k_1)\sn(s_1^{12},k_4)|^2|\cos s_2^3\cos s_2^8|^2+ \\ 
& + |\dn(s_3^7,k_3)||\dn(s_3^{12},k_6)|^2|\cos s_2^3\cos s_2^8|^2, 
\end{split}
\ee  

\be \label{3.8} 
\begin{split}
& (2.9)\times (2.13) \Rightarrow \\ 
& \frac{|s_1^4|^{2n_1}|s_1^9|^{2n_4}}{|\Gamma(s_2^5)\Gamma(s_2^{10})|^2}+ 
\frac{|s_3^4|^{n_3}|s_3^9|^{2n_5}}{|\Gamma(s_2^5)\Gamma(s_2^{10})|^2}= \\ 
& = \frac{|s_2^4|^{2n_2}|s_2^9|^{2n_6}}{|\Gamma(s_1^5)\Gamma(s_1^{10})|^2}+ 
\frac{|s_2^4|^{2n_2}|s_2^9|^{2n_6}}{|\Gamma(s_3^5)||\Gamma(s_3^{10})|^2}, 
\end{split}
\ee  

\be \label{3.9} 
\begin{split}
& (2.9)\times (2.18) \Rightarrow \\  
& |J_{p_2}(s_2^6)J_{p_5}(s_2^{11})|^2|s_1^4|^{2n_1}|s_1^9|^{2n_4}+ \\ 
& + |J_{p_2}(s_2^6)J_{p_5}(s_2^{11})|^2|s_3^4|^{n_3}|s_3^9|^{2n_5}= \\ 
& = |J_{p_1}(s_1^6)J_{p_4}(s_1^{11})|^2|s_2^4|^{2n_2}|s_2^9|^{2n_6}+ \\ 
& + |J_{p_3}(s_3^6)||J_{p_6}(s_3^{11})|^2|s_2^4|^{2n_2}|s_2^9|^{2n_6}, 
\end{split}
\ee  

\be \label{3.10} 
\begin{split}
& (2.9)\times (2.23) \Rightarrow \\  
& |\cn(s_2^7,k_2)\cn(s_2^{12},k_5)|^2|s_1^4|^{2n_1}|s_1^9|^{2n_4}+ \\ 
& + |\cn(s_2^7,k_2)\cn(s_2^{12},k_5)|^2|s_3^4|^{n_3}|s_3^9|^{2n_5}= \\ 
& = |\sn(s_1^7,k_1)\sn(s_1^{12},k_4)|^2|s_2^4|^{2n_2}|s_2^9|^{2n_6}+ \\ 
& + |\dn(s_3^7,k_3)||\dn(s_3^{12},k_6)|^2|s_2^4|^{2n_2}|s_2^9|^{2n_6}, 
\end{split}
\ee  

\be \label{3.11} 
\begin{split}
& (2.13)\times (2.18) \Rightarrow \\ 
& \frac{|J_{p_2}(s_2^6)J_{p_5}(s_2^{11})|^2}{|\Gamma(s_1^5)\Gamma(s_1^{10})|^2}+ 
\frac{|J_{p_2}(s_2^6)J_{p_5}(s_2^{11})|^2}{|\Gamma(s_3^5)||\Gamma(s_3^{10})|^2}= \\ 
& = \frac{|J_{p_1}(s_1^6)J_{p_4}(s_1^{11})|^2}{|\Gamma(s_2^5)\Gamma(s_2^{10})|^2}+ 
\frac{|J_{p_3}(s_3^6)||J_{p_6}(s_3^{11})|^2}{|\Gamma(s_2^5)\Gamma(s_2^{10})|^2}, 
\end{split}
\ee  

\be \label{3.12} 
\begin{split}
& (2.13)\times (2.23) \Rightarrow \\  
& \frac{|\cn(s_2^7,k_2)\cn(s_2^{12},k_5)|^2}{|\Gamma(s_1^5)\Gamma(s_1^{10})|^2}+ 
\frac{|\cn(s_2^7,k_2)\cn(s_2^{12},k_5)|^2}{|\Gamma(s_3^5)||\Gamma(s_3^{10})|^2}= \\ 
& = \frac{|\sn(s_1^7,k_1)\sn(s_1^{12},k_4)|^2}{|\Gamma(s_2^5)\Gamma(s_2^{10})|^2}+ 
\frac{|\dn(s_3^7,k_3)||\dn(s_3^{12},k_6)|^2}{|\Gamma(s_2^5)\Gamma(s_2^{10})|^2}, 
\end{split}
\ee 

\be \label{3.13} 
\begin{split}
& (2.18)\times (2.23) \Rightarrow \\ 
& |\cn(s_2^7,k_2)\cn(s_2^{12},k_5)|^2|J_{p_1}(s_1^6)J_{p_4}(s_1^{11})|^2+ \\
& + |\cn(s_2^7,k_2)\cn(s_2^{12},k_5)|^2|J_{p_3}(s_3^6)||J_{p_6}(s_3^{11})|^2= \\ 
& = |\sn(s_1^7,k_1)\sn(s_1^{12},k_4)|^2|J_{p_2}(s_2^6)J_{p_5}(s_2^{11})|^2+ \\ 
& + |\dn(s_3^7,k_3)||\dn(s_3^{12},k_6)|^2|J_{p_2}(s_2^6)J_{p_5}(s_2^{11})|^2. 
\end{split}
\ee 
\end{mydef1} 

\begin{remark}
The following is true: the set of level curves 
\bdis 
\overset{(n)}{\Omega}_l,\ n=8,\dots,12,\ l=1,2,3 
\edis  
is defined as the set 
\bdis 
\{(2.2), (2.7), (2.11), (2.16), (2.21)\}
\edis  
of loci by means of the values 
\bdis 
c_l=\left|\zeta\left(\frac 12+i\alpha_1^{l,1}\right)\right|,\ l=1,2,3 
\edis  
(comp. \cite{8}, (3.1), (3.2)), where (see \cite{9}, (1.10)) 
\bdis 
\alpha_1^{l,1}=\alpha_1(U,\pi L,f_l;[|\zeta_{0,5}|^2]). 
\edis  
Since 
\bdis 
c_l=\left|\zeta\left(\frac 12+i\alpha_1^{l,1}\right)\right|=
\left|\zeta\left(\frac 12+i\alpha_1(U,\pi L,f_l;[|\zeta_{0,5}|^2])\right)\right|, 
\edis 
then we have $c_l$'s -- as initial conditions --  twice determined by the function 
\bdis 
\left|\zf\right|. 
\edis 
\end{remark}

I would like to thank Michal Demetrian for his moral support of my study of Jacob's ladders.


\begin{thebibliography}{29}% Replace 9 by 99 if 10 or more references
\bibitem{1}
J. Moser,
`Jacob's ladders and almost exact asymptotic representation of the Hardy-Littlewood integral`,
Math. Notes 88, (2010), 414-422, arXiv: 0901.3937.
%
\bibitem{2}
J. Moser,
`Jacob's ladders, structure of the Hardy-Littlewood integral and some new class of nonlinear integral equations`,
Proc. Steklov Inst. 276 (2011), 208-221, arXiv: 1103.0359.
%
\bibitem{3}
J. Moser, `Jacob's ladders, reverse iterations and new infinite set of $L_2$-orthogonal systems generated by the
Riemann zeta-function, arXiv: 1402.2098.
%
\bibitem{4}
J. Moser, `Jacob's ladders, factorization and metamorphoses as an appendix to the Riemann functional equation for $\zeta(s)$ on
the critical line`, Proc. Steklov Inst. 296 (2017), pp. 92-102, arXiv: 1506.00442v1. 
%
\bibitem{5}
J. Moser, 'Jacob's ladders, interactions between $\zeta$-oscillating systems and $\zeta$-analogue of an elementary trigonometric identity', arXiv: 1609.09293v1, Proc. Steklov Inst. 299, 189-204, 2017. 
%
\bibitem{6}
J. Moser, `Jacob ladders, crossbreeding, secondary crossbreeding and synergetic phenomena generated by the Riemann's zeta-function and some elementary functions on disconnected sets of the critical line`, arXiv: 1806.07095v1. 
\bibitem{7}
J. Moser, `Jacob's ladders and grafting of the complete hybrid formulas into $\zeta$-synergetic meta-functional equation for the Riemann's zeta-function`, arXiv: 1809.05327v1. 
\bibitem{8}
J. Moser, `Jacob ladders and infinite set of transmutations of asymptotic complete hybrid formula on level curves in Gauss' plane`, arXiv: 1905.06078v1. 
\bibitem{9}
J. Moser, `Jacob's ladders and exact meta-functional equations on level curves as global quantitative characteristics of synergetic phenomenons excited by the functions $|\zf|^2$`, arXiv: 1906.02440v1. 


\end{thebibliography}
\end{document}